\documentclass[11pt,reqno]{amsart}

\usepackage{amsmath}
\usepackage{amssymb}
\usepackage{amsthm}

\newtheorem{theorem}{Theorem}[section]
\newtheorem{prop}[theorem]{Proposition}
\newtheorem{lem}[theorem]{Lemma}
\newtheorem*{cor}{Corollary}
\theoremstyle{definition}
\newtheorem{defn}[theorem]{Definition}
\theoremstyle{remark}
\newtheorem*{rem}{Remark}
\numberwithin{equation}{section}

\begin{document}

\title[Fourier transform on locally compact quantum groups]
{Fourier transform on locally compact quantum groups}

\author{Byung-Jay Kahng}
\date{}
\address{Department of Mathematics and Statistics\\ Canisius College\\
Buffalo, NY 14208}
\email{kahngb@canisius.edu}

\begin{abstract}
The notion of Fourier transform is among the more important tools in analysis, which 
has been generalized in abstract harmonic analysis to the level of abelian locally 
compact groups.  The aim of this paper is to further generalize the Fourier transform: 
Motivated by some recent works by Van Daele in the multiplier Hopf algebra framework, 
and by using the Haar weights, we define here the (generalized) Fourier transform and 
the inverse Fourier transform, at the level of locally compact quantum groups.  We will 
then consider the analogues of the Fourier inversion theorem, Plancherel theorem, and 
the convolution product.  Along the way, we also obtain an alternative description of 
the dual pairing map between a quantum group and its dual.
\end{abstract}
\maketitle

\section{Introduction}

The Fourier transform has been known for quite some time, and is among the very powerful 
tools in classical analysis.  In abstract harmonic analysis (see \cite{HR}, \cite{Fo}), 
the theory was generalized to the case of abelian locally compact groups.

Let us briefly summarize: Given a locally compact abelian (LCA) group $G$, its dual
object $\hat{G}$ is the group of continuous, $\mathbb{T}$-valued characters on $G$. 
The dual group $\hat{G}$ can be given an appropriate topology, making it into an LCA 
group again.  By Pontryagin duality, it is known that the dual of $\hat{G}$ is isomorphic 
to $G$.  Using the Haar measure on $G$, it is then possible to define the Fourier 
transform of a continuous function having compact support $f\in C_c(G)$ (or, even 
a Schwartz function), obtaining $\hat{f}\in C_0(\hat{G})$.  Of course, one can begin 
with the functions on $\hat{G}$, and define the inverse Fourier transform.  The Fourier 
inversion theorem holds, as well as Plancherel's theorem.  One by-product is that we 
obtain in this way a Hilbert space isomorphism $L^2(G)\cong L^2(\hat{G})$.  Refer to 
the standard textbooks on the theory, for instance, \cite{Fo}.

If the group $G$ is non-abelian, we no longer can define $\hat{G}$ as the dual group.
Tannaka, Krein, and others have been able to define a dual object to a non-abelian
group, from which the original group can be recovered.  However, with the dual object 
not being a group, the Pontryagin duality does not hold anymore, and it is not possible 
to define the Fourier transform and the inverse Fourier transform between $G$ and 
$\hat{G}$.  A modified version of the Fourier transform does exist, and is being used, 
in the representation theory of non-abelian, compact groups.  But still, from the 
duality point of view, it is not really satisfactory.

The environment is now better, with the recent development of the theory of locally 
compact quantum groups \cite{KuVa}, \cite{KuVavN}, \cite{MNW}, \cite{VDvN}.  We now 
know that Pontryagin duality can be naturally extended to the wider setting of locally 
compact quantum groups.  Indeed, given a locally compact quantum group $(M,\Delta)$, 
its dual object $(\hat{M},\hat{\Delta})$ is also a locally compact quantum group, and 
moreover, the dual of $\hat{M}$ is isomorphic to $M$.

On the other hand, a workable notion of Fourier transform at the quantum group level 
has been lacking so far.  Among the challenges is that because of the way the dual 
quantum group is constructed, the Hilbert spaces for $M$ and $\hat{M}$ are identical, 
and the Fourier transform tends to be ``hidden'' (It is essentially like identifying 
the spaces $L^2(G)=L^2(\hat{G})$, making the Fourier transform irrelevant.).

In his recent works, Van Daele has been working to improve this situation.  In 
particular, in his preprint \cite{VDFourier} (see also \cite{VDaqg}), he proposes 
a definition for the generalized Fourier transform in the setting of multiplier Hopf 
algebras and algebraic quantum groups.  He suggests in the paper that more can be 
done at the operator algebra level of locally compact quantum groups, and indeed, 
some indications of these new developments do appear in his other recent papers 
(for instance, \cite{VDvN}).  However, at the time of writing the present paper, 
the author does not know of any place where his results are actually written out.

The present paper grew out in the hope of filling this gap.  Using Van Daele's 
approach in \cite{VDFourier} as a guideline, and by taking advantage of the 
generalized Pontryagin duality, we define and explore here the notion of Fourier 
transform in the setting of locally compact quantum groups.  The point we wish to 
make is that even though the Fourier transform may seem hidden, it is still there, 
and we are proposing a way to harness its usefulness.  This generalized Fourier 
transform will be shown to satisfy many of the familiar results from classical 
analysis.

Due to the technical differences between the setting of multiplier Hopf algebras 
and the framework of locally compact quantum groups, the proofs and the details 
are quite different, while the end-results may look similar.  Nevertheless, the 
author can never downplay the motivation and the strong influence Van Daele's 
works gave him in preparing the present paper.  

Here is how this paper is organized.  We begin, in Section 2, by recalling the
definition and some main results on locally compact quantum groups.  We chose 
to work with the von Neumann algebra approach (as in \cite{KuVavN}), which is 
known to be equivalent to the $C^*$-algebraic framework.  Special attentions are 
given to Haar weights, the multiplicative unitary operator, and the antipode maps.

In Section 3, we give the definition of the generalized Fourier transform, sending 
certain elements of $M$ to elements in $\hat{M}$.  We will also define the inverse 
Fourier transform, and observe the analogues of the Fourier inversion theorem 
and the Plancherel theorem.  A formula for the ``convolution product'' will be 
also obtained.

In Section 4, we will use the Fourier transform to give an alternative (and 
apparently new) description of the dual pairing between certain dense subalgebras 
of $M$ and $\hat{M}$.  This will be useful in our future works.  Finally, we 
added a brief Appendix (Section~5), where we observe how all this is reflected 
in the special case of an ordinary locally compact group.

\section{Preliminaries: von Neumann algebraic quantum groups and Haar weights}

We will use the standard notations from the theory of weights (see \cite{St}, 
\cite{Tk2}).  The weights we will be working with are normal, semi-finite faithful 
weights (``n.s.f. weights'', for short) on von Neumann algebras.  For an n.s.f. 
weight $\varphi$ on a von Neumann algebra $M$, we write:
\begin{itemize}
  \item ${\mathfrak M}_{\varphi}^+=\bigl\{x\in M^+:\varphi(x)<\infty\bigr\}$
  \item ${\mathfrak N}_{\varphi}=\{x\in M:x^*x\in{\mathfrak M}_{\varphi}^+\}$
  \item ${\mathfrak M}_{\varphi}=\{\sum_{i=1}^{n}y_i^*x_i:x_1,\dots,x_n,y_1,\dots,y_n
\in{\mathfrak N}_{\varphi}\}$
\end{itemize}
The space ${\mathfrak M}_{\varphi}$ is a ${}^*$-subalgebra of $M$, which is the 
``definition domain'' of the weight $\varphi$.  Observe that ${\mathfrak M}_{\varphi}$
is obtained as the linear span of ${\mathfrak M}_{\varphi}^+$ in $M$.

Let us begin with the definition of a von Neumann algebraic locally compact quantum 
group, as given by Kustermas and Vaes \cite{KuVavN}.  This definition is known to 
be equivalent to the definition in the $C^*$-algebra setting \cite{KuVa}, \cite{MNW}.  
As in the case of the $C^*$-algebraic quantum groups, the existence of Haar (invariant) 
weights is assumed as a part of the definition.  The noticeable difference between 
the two approaches is the absence of the density conditions in the von Neumann 
algebra setting: It turns out that they follow automatically from the other 
conditions.  Refer also to the recent paper by Van Daele \cite{VDvN}, which gives 
an improved approach to the subject and is more natural.

\begin{defn}\label{lcqgvN}
Let $M$ be a von Neumann algebra, together with a unital normal ${}^*$-homomorphism
$\Delta:M\to M\otimes M$ such that the ``coassociativity condition'' holds: 
$(\Delta\otimes\operatorname{id})\Delta=(\operatorname{id}\otimes\Delta)\Delta$.
Furthermore, we assume the existence of a left invariant weight and a right invariant
weight, as follows:
\begin{itemize}
  \item $\varphi$ is an n.s.f. weight on $M$ that is left invariant:
$$
\varphi\bigl((\omega\otimes\operatorname{id})(\Delta x)\bigr)=\varphi(x)\omega(1),
\quad {\text {for all }}\omega\in M_*^+,\  x\in{\mathfrak M}_{\varphi}^+.
$$
  \item $\psi$ is an n.s.f. weight on $M$ that is right invariant:
$$
\psi\bigl((\operatorname{id}\otimes\omega)(\Delta x)\bigr)=\psi(x)\omega(1),
\quad {\text {for all }}\omega\in M_*^+,\  x\in{\mathfrak M}_{\psi}^+.
$$
\end{itemize}
Then we call $(M,\Delta)$ a {\em von Neumann algebraic quantum group\/}.
It can be shown that the Haar weights are unique, up to scalar multiplication.
\end{defn}

Let us fix $\varphi$.  By means of the GNS-construction $({\mathcal H},\iota,\Lambda)$ 
for $\varphi$, we view $M$ as a subalgebra of the operator algebra ${\mathcal B}
({\mathcal H})$, such as $M=\iota(M)\subseteq{\mathcal B}({\mathcal H})$.  
So we will have: $\bigl\langle\Lambda(x),\Lambda(y)\bigr\rangle=\varphi(y^*x)$ 
for $x,y\in{\mathfrak N}_{\varphi}$, and $x\Lambda(y)=\Lambda(xy)$ for 
$y\in{\mathfrak N}_{\varphi}$, $x\in M$.  As in standard weight theory, we can also
consider the modular conjugation and the modular automorphism group of $\varphi$.

Meanwhile, there exists a unitary operator $W\in{\mathcal B}({\mathcal H}\otimes
{\mathcal H})$, called the {\em multiplicative unitary operator\/} for $(M,\Delta)$.
The operator $W$ is defined by $W^*\bigl(\Lambda(x)\otimes\Lambda(y)\bigr)
=(\Lambda\otimes\Lambda)\bigl((\Delta y)(x\otimes1)\bigr)$, for $x,y\in
{\mathfrak N}_{\varphi}$.  It satisfies the pentagon equation of Baaj and Skandalis
\cite{BS}: $W_{12}W_{13}W_{23}=W_{23}W_{12}$, and one can check that $\Delta x
=W^*(1\otimes x)W$, for all $x\in M$.  It is essentially the ``left regular 
representation'' (associated with $\varphi$), and it also gives us the following 
useful characterization of $M$:
\begin{equation}\label{(M)}
M=\overline{\bigl\{(\operatorname{id}\otimes\omega)(W):\omega\in{\mathcal B}
({\mathcal H})_*\bigr\}}^w\,\bigl(\subseteq{\mathcal B}
({\mathcal H})\bigr),
\end{equation}
where ${-}^w$ denotes the von Neumann algebra closure (for instance, the closure
under $\sigma$-weak topology).

There are other possible (and useful) characterizations of $(M,\Delta)$.  See 
\cite{KuVavN}.  Meanwhile, if we take the norm closure in equation \eqref{(M)}, 
instead of the weak closure, we obtain a $C^*$-algebra $A$.  It turns out that 
by restricting $\Delta$ to $A$, we obtain a reduced $C^*$-algebraic quantum group 
$(A,\Delta)$.  Going the other way, we could begin with a $C^*$-algebraic quantum 
group and obtain a von Neumann algebraic quantum group by taking the weak closure 
of the underlying $C^*$-algebra in the GNS Hilbert space of a left Haar weight.

The main point in all these is that the above Definition \ref{lcqgvN} is a valid
definition of a locally compact quantum group.  For instance, the existence of 
other quantum group structure maps like ``antipode'' can be proved from the 
defining axioms.

Constructing the antipode is quite technical (it uses the right Haar weight), 
and we refer the reader to the main papers \cite{KuVa}, \cite{KuVavN}.  See 
also an improved treatment given in \cite{VDvN}, where the antipode is defined
in a more natural way by means of Tomita--Takesaki theory.  For our purposes, 
we will work with the following useful characterization of the antipode map $S$:
\begin{equation}\label{(antipode)}
S\bigl((\operatorname{id}\otimes\omega)(W)\bigr)
=(\operatorname{id}\otimes\omega)(W^*).
\end{equation}
In fact, the subspace consisting of the elements $(\operatorname{id}\otimes\omega)
(W)$, for $\omega\in{\mathcal B}({\mathcal H})_*$, is dense in $M$, and forms 
a core for $S$.  Meanwhile, there exist a unique ${}^*$-antiautomorphism $R$ 
(called the ``unitary antipode'') and a unique continuous one parameter group 
$\tau$ on $M$ (called the ``scaling group'') such that we have: 
$S=R\tau_{-\frac{i}{2}}$.

Since $(R\otimes R)\Delta=\Delta^{\operatorname{cop}}R$, where 
$\Delta^{\operatorname{cop}}$ is the co-opposite comultiplication (i.\,e. 
$\Delta^{\operatorname{cop}}=\chi\circ\Delta$, for $\chi$ the flip map 
on $M\otimes M$), the weight $\varphi\circ R$ is right invariant.  So for 
convenience, we will from now on choose $\psi$ to equal $\varphi\circ R$. 
The GNS map for $\psi$ will be written as $\Gamma$.

Next, let us consider the dual quantum group.  Working with the other leg of 
the multiplicative unitary operator $W$ than in equation \eqref{(M)}, we define:
\begin{equation}\label{(Mhat)}
\hat{M}=\overline{\bigl\{(\omega\otimes\operatorname{id})(W):\omega
\in{\mathcal B}({\mathcal H})_*\bigr\}}^w\,\bigl(\subseteq{\mathcal B}
({\mathcal H})\bigr).
\end{equation}
This is indeed a von Neumann algebra.  The comultiplication on it is defined by
$\hat{\Delta}(y)=\Sigma W(y\otimes1)W^*\Sigma$, for all $y\in\hat{M}$ (here, 
$\Sigma$ is the flip on ${\mathcal H}\otimes{\mathcal H}$).  The general theory 
assures that $(\hat{M},\hat{\Delta})$ is again a von Neumann algebraic quantum 
group, together with appropriate Haar weights $\hat{\varphi}$ and $\hat{\psi}$. 
The operator $\Sigma W^*\Sigma$ is the multiplicative unitary for $(\hat{M},
\hat{\Delta})$.  It  can be shown that $W\in M\otimes\hat{M}$ and $\Sigma W^*
\Sigma\in\hat{M}\otimes M$.

The left Haar weight $\hat{\varphi}$ on $(\hat{M},\hat{\Delta})$ is uniquely 
characterized by the GNS data $({\mathcal H},\iota,\hat{\Lambda})$, where the 
GNS map $\hat{\Lambda}:{\mathfrak N}_{\hat{\varphi}}\to{\mathcal H}$ is given 
by the following formulas (See Proposition~8.14 of \cite{KuVa}):
\begin{equation}\label{(phihat)}
\hat{\Lambda}\bigl((\omega\otimes\operatorname{id})(W)\bigr)=\xi(\omega)\quad
{\text {and}}\quad\bigl\langle\xi(\omega),\Lambda(x)\bigr\rangle=\omega(x^*).
\end{equation}
To be a little more precise, consider:
$$
{\mathcal I}=\bigl\{\omega\in{\mathcal B}({\mathcal H})_*:\exists L\ge0\ 
{\text {such that }}|\omega(x^*)|\le L\|\Lambda(x)\|\ {\text {for all }}
x\in{\mathfrak N}_{\varphi}\bigr\}.
$$
Then for every $\omega\in{\mathcal I}$, we can find $\xi(\omega)\in{\mathcal H}$
such that $\omega(x^*)=\bigl\langle\xi(\omega),\Lambda(x)\bigr\rangle$ for all 
$x\in{\mathfrak N}_{\varphi}$ (by Riesz theorem).  The equation~\eqref{(phihat)} 
above is understood as saying that the elements $(\omega\otimes\operatorname
{id})(W)$, $\omega\in{\mathcal I}$, form a core for $\hat{\Lambda}$ and that 
$\hat{\Lambda}\bigl((\omega\otimes\operatorname{id})(W)\bigr)=\xi(\omega)$.
See also \cite{VDvN}, to learn more on this construction.

Meanwhile, analogously as in equation~\eqref{(antipode)}, with $\Sigma W^*
\Sigma$ now being the multiplicative unitary, the (dense) subspace of the 
elements $(\omega\otimes\operatorname{id})(W^*)$, for $\omega\in{\mathcal B}
({\mathcal H})_*$, forms a core for the antipode $\hat{S}$, and $\hat{S}$ is 
characterized by 
\begin{equation}\label{(antipodehat)}
\hat{S}\bigl((\omega\otimes\operatorname{id})(W^*)\bigr)
=(\omega\otimes\operatorname{id})(W).
\end{equation}
The unitary antipode and the scaling group can be also found, giving us the
polar decomposition of $\hat{S}=\hat{R}\hat{\tau}_{-\frac{i}2}$.  As before,
we may fix the right Haar weight as $\hat{\psi}=\hat{\varphi}\circ\hat{R}$, 
with the corresponding GNS map written as $\hat{\Gamma}$.

Repeating the whole process beginning with $({\mathcal H},\iota,\hat{\Lambda})$, 
we can further construct the dual $(\hat{\hat{M}},\hat{\hat{\Delta}})$ of 
$(\hat{M},\hat{\Delta})$.  The generalized Pontryagin duality result (see 
\cite{KuVa}, \cite{KuVavN}) says: $(\hat{\hat{M}},\hat{\hat{\Delta}})=
(M,\Delta)$, with $\hat{\hat{\varphi}}=\varphi$ and $\hat{\hat{\Lambda}}
=\Lambda$.  The other structure maps for $(\hat{\hat{M}},\hat{\hat{\Delta}})$ 
are also identified with those of $(M,\Delta)$: For instance, $\hat{\hat{S}}=S$, 
$\hat{\hat{\psi}}=\psi$, etc.  One useful result is the following, similar to 
equation \eqref{(phihat)} above (See Proposition~8.30 of \cite{KuVa}, with 
$\pi$ now read as the embedding map 
$\iota$ and $\hat{\hat{\Lambda}}=\Lambda$):
\begin{equation}\label{(phihatdual)}
\bigl\langle\Lambda\bigl((\operatorname{id}\otimes\omega)(W^*)\bigr),
\hat{\Lambda}(y)\bigr\rangle=\omega(y^*).
\end{equation}
Here again, we actually need to consider a set $\hat{\mathcal I}$ similarly
defined as before, and the equation \eqref{(phihatdual)} is accepted with 
the understanding that the elements $(\operatorname{id}\otimes\omega)(W^*)$, 
$\omega\in\hat{\mathcal I}$, form a core for $\Lambda$.

We wrap up the section here.  For further details, we refer the reader to 
the fundamental papers on the subject: \cite{BS}, \cite{Wr7}, \cite{KuVa}, 
\cite{KuVavN}, \cite{MNW}, \cite{VDvN}.

\section{The generalized Fourier transform}

Let us denote by ${\mathcal A}$ and $\hat{\mathcal A}$, the dense subalgebras 
of $M$ and $\hat{M}$, similar to the ones given in equations \eqref{(M)} and 
\eqref{(Mhat)}:
\begin{align}
{\mathcal A}&=\bigl\{(\operatorname{id}\otimes\omega)(W):\omega\in \hat{M}_*
\bigr\}\,(\subseteq M),  \notag \\
\hat{\mathcal A}&=\bigl\{(\omega\otimes\operatorname{id})(W):\omega\in M_*
\bigr\}\,(\subseteq\hat{M}).  \notag
\end{align}
The fact that these are indeed subalgebras follows from the fact that $W\in M
\otimes\hat{M}$ and that $W$ is a multiplicative unitary operator (see \cite{BS}, 
\cite{Wr7}).  See the lemma below:

\begin{lem}\label{lemsubalg}
By the multiplicativity of $W$ (the pentagon equation), the following results hold:
\begin{enumerate}
 \item For $\omega_1,\omega_2\in M_*$, we have: $(\omega_1\otimes
\operatorname{id})(W)(\omega_2\otimes\operatorname{id})(W)
=(\mu\otimes\operatorname{id})(W)$, where $\mu\in M_*$ is such that
$\mu(x)=(\omega_1\otimes\omega_2)(\Delta x)$, for $x\in M$.
 \item For $\theta_1,\theta_2\in\hat{M}_*$, we have: $(\operatorname{id}
\otimes\theta_1)(W)(\operatorname{id}\otimes\theta_2)(W)
=(\operatorname{id}\otimes\nu)(W)$, where $\nu\in\hat{M}_*$ is such that
$\nu(y)=(\theta_1\otimes\theta_2)\bigl(\hat{\Delta}^{\operatorname{cop}}(y)
\bigr)$, for $y\in\hat{M}$.
\item For $\theta_1,\theta_2\in\hat{M}_*$, we have: $(\operatorname{id}
\otimes\theta_1)(W^*)(\operatorname{id}\otimes\theta_2)(W^*)
=(\operatorname{id}\otimes\nu')(W^*)$, where $\nu'\in\hat{M}_*$ is such that
$\nu'(y)=(\theta_1\otimes\theta_2)\bigl(\hat{\Delta}(y)\bigr)$, for $y\in\hat{M}$.
\end{enumerate}
\end{lem}

\begin{rem}
By (1), we see that $\hat{\mathcal A}$ is a subalgebra of $\hat{M}$, while (2) 
shows that ${\mathcal A}$ is a subalgebra of $M$.  (3) will be also useful later.
\end{rem}

\begin{proof}
By the pentagon equation ($W_{12}W_{13}W_{23}=W_{23}W_{12}$), 
we have:
\begin{align}
(\omega_1\otimes\operatorname{id})(W)(\omega_2\otimes\operatorname{id})(W)
&=(\omega_1\otimes\omega_2\otimes\operatorname{id})(W_{13}W_{23})
\notag \\
&=(\omega_1\otimes\omega_2\otimes\operatorname{id})(W_{12}^*W_{23}W_{12})
=(\mu\otimes\operatorname{id})(W),
\notag
\end{align}
where $\mu\in M_*$ is such that $\mu(x)=(\omega_1\otimes\omega_2)
\bigl(W^*(1\otimes x)W\bigr)$.  Remembering the definition of $\Delta x
=W^*(1\otimes x)W$, we obtain the result (1).  See \cite{BS} for the same
result.  The proof for the other two are similar. 
\end{proof}

As for the involution, the following lemma will be useful (see \cite{KuVavN}):

\begin{lem}\label{lemmasharp}
We will write $\omega\in M_*^{\sharp}$, if $\omega\in M_*$ is such that there 
exists an element $\omega^{\sharp}\in M_*$, given by: $\omega^{\sharp}(x)
=\bar{\omega}\bigl(S(x)\bigr)=\overline{\omega\bigl([S(x)]^*\bigr)}$, 
for all $x\in{\mathcal D}(S)$.  In that case, we have:
$$
\bigl((\omega\otimes\operatorname{id})(W)\bigr)^*
=(\omega^{\sharp}\otimes\operatorname{id})(W).
$$
The choice of $\omega^{\sharp}$ above is unique, in the sense that if there 
exists a $\rho\in M_*$ satisfying $\bigl((\omega\otimes\operatorname{id})(W)
\bigr)^*=(\rho\otimes\operatorname{id})(W)$, then we have: $\rho=\omega^{\sharp}$.

Meanwhile, the subspace $M_*^{\sharp}\,(\subseteq M_*)$ is a dense subalgebra 
of $M_*$ (in the sense of (2) of Lemma~\ref{lemsubalg}), and is closed under 
taking ${}^{\sharp}$.
\end{lem}

\begin{proof}
If $\omega\in M_*^{\sharp}$, then by the characterization of the antipode $S$
given in equation \eqref{(antipode)}, we have, for any $\theta\in{\mathcal B}
({\mathcal H})_*$:
\begin{align}
&\theta\bigl((\omega^{\sharp}\otimes\operatorname{id})(W)\bigr)
=\omega^{\sharp}\bigl((\operatorname{id}\otimes\theta)(W)\bigr)
=\overline{\omega\bigl([S((\operatorname{id}\otimes\theta)(W))]^*\bigr)}
\notag \\
&\quad=\overline{\omega\bigl([(\operatorname{id}\otimes\theta)(W^*)]^*\bigr)}
=\overline{\omega\bigl((\operatorname{id}\otimes\bar{\theta})(W)\bigr)}
=\overline{\bar{\theta}\bigl((\omega\otimes\operatorname{id})(W)\bigr)}
\notag \\
&\quad=\theta\bigl([(\omega\otimes\operatorname{id})(W)]^*).
\notag
\end{align}
Since $\theta$ is arbitrary, we see that 
$(\omega^{\sharp}\otimes\operatorname{id})(W)
=\bigl[(\omega\otimes\operatorname{id})(W)\bigr]^*$.

For uniqueness, suppose $\bigl((\omega\otimes\operatorname{id})(W)\bigr)^*
=(\rho\otimes\operatorname{id})(W)$.  Then for any $x=(\operatorname{id}
\otimes\theta)(W)$, $\theta\in{\mathcal B}({\mathcal H})_*$, we have:
\begin{align}
\rho(x)&=\theta\bigl((\rho\otimes\operatorname{id})(W)\bigr)
=\theta\bigl([(\omega\otimes\operatorname{id})(W)]^*\bigr)
=\overline{\bar{\theta}\bigl((\omega\otimes\operatorname{id})(W)\bigr)}
\notag \\
&=\overline{\omega\bigl((\operatorname{id}\otimes\bar{\theta})(W)\bigr)}
=\overline{\omega\bigl([(\operatorname{id}\otimes\theta)(W^*)]^*\bigr)}
=\overline{\omega\bigl([S(x)]^*\bigr)}=\omega^{\sharp}(x).
\notag
\end{align}
Since the elements $(\operatorname{id}\otimes\theta)(W)$, $\theta\in{\mathcal B}
({\mathcal H})_*$, form a core for ${\mathcal D}(S)$, the result holds for 
all $x\in{\mathcal D}(S)$.

Finally, it is easy to see that $M_*^{\sharp}$ is a subalgebra of $M_*$, 
and is closed under taking ${}^{\sharp}$ (note that $(\omega^{\sharp})^{\sharp}
=\omega$).  To see if $M_*^{\sharp}$ is dense in $M_*$, note that for any 
analytic element $\omega\in M_*$, we can find $\omega^{\sharp}$ by 
$\omega^{\sharp}=\bar{\omega}\circ{\tau}_{-\frac{i}2}\circ R\,(=\bar{\omega}
\circ S)$.  Since such analytic elements are dense, we have our proof.
\end{proof}

\begin{rem}
For convenience, we will from now on use the notation $\lambda(\omega)$, 
$\omega\in M_*$, to mean $\lambda(\omega)=(\omega\otimes
\operatorname{id})(W)\,\in\hat{\mathcal A}\,(\subseteq\hat{M})$.  So the 
last part of Lemma~\ref{lemmasharp} implies that $\lambda(M_*^{\sharp})$ 
is dense in $\hat{\mathcal A}$ (and $\hat{M}$).  On the other hand, since 
$S$ may be unbounded in general, it may not be all of $\hat{\mathcal A}$.

Dually, let us introduce also the notation $\hat{\lambda}(\theta)$, $\theta
\in\hat{M}_*$, to mean $\hat{\lambda}(\theta)=(\theta\otimes\operatorname{id})
(\Sigma W^*\Sigma)=(\operatorname{id}\otimes\theta)(W^*)\,\in M$.  We can 
define $\hat{M}_*^{\sharp}$ the same way as above: i.\,e. $\theta\in
\hat{M}_*^{\sharp}$ if $\bigl((\operatorname{id}\otimes\theta)(W^*)\bigr)^*
=(\operatorname{id}\otimes\theta^{\sharp})(W^*)$, for $\theta^{\sharp}
=\bar{\theta}\circ\hat{S}$.  Since we have $(\operatorname{id}\otimes
\theta^{\sharp})(W^*)=(\operatorname{id}\otimes\bar{\theta})(W)$, we see that 
$\hat{\lambda}(\hat{M}_*^{\sharp})$ is contained and dense in ${\mathcal A}$ 
(and $M$).
\end{rem}

Since we will be formulating the Fourier transform in terms of the Haar weights,
let us recall the definition of the spaces ${\mathcal I}$ and $\hat{\mathcal I}$,
as mentioned in the previous section (see also \cite{KuVavN}):
\begin{align}
{\mathcal I}&=\bigl\{\omega\in M_*:\exists L\ge0 {\text { such that }}
|\omega(x^*)|\le L\|\Lambda(x)\|, \forall x\in{\mathfrak N}_{\varphi}\bigr\}, 
\notag \\
\hat{\mathcal I}&=\bigl\{\theta\in\hat{M}_*:\exists L\ge0 {\text { such that }}
|\theta(y^*)|\le L\|\hat{\Lambda}(y)\|, 
\forall y\in{\mathfrak N}_{\hat{\varphi}}\bigr\}.  \notag
\end{align}
As noted in the previous section, the space $\lambda({\mathcal I})\,(\subseteq
\hat{M})$ is a core for $\hat{\Lambda}$, while the space $\hat{\lambda}
(\hat{\mathcal I})\,(\subseteq M)$ is a core for $\Lambda$.

Let us now turn our attention to our main goal of defining the Fourier transform. 
In \cite{VDaqg} and \cite{VDFourier}, at the level of multiplier Hopf algebras, 
the Fourier transform of an element $a$ is defined as the linear functional 
$\omega=\varphi(\cdot\,a)$.  In fact, the dual multiplier Hopf algebra is 
characterized as the collection of such linear functionals.  With this as 
motivation, and considering that the multiplicative unitary operator $W$ 
provides the duality (see comments in Section 3 of \cite{VDvN}, and see also 
equation~\eqref{(pairing)} below), we propose to take the definition of the 
Fourier transform as ${\mathcal F}(a)=(\omega\otimes\operatorname{id})(W)$, 
where $\omega=\varphi(\cdot\,a)$.  See the formulation given in Definition~\ref
{FTdefn} below.

\begin{defn}\label{FTdefn}
For $a\in\hat{\lambda}(\hat{\mathcal I})\,(\subseteq M)$, define ${\mathcal F}(a)
\in\hat{M}$, such that
$$
{\mathcal F}(a):=(\varphi\otimes\operatorname{id})\bigl(W(a\otimes1)\bigr).
$$
Note that formally, we can write it as ${\mathcal F}(a)=(\omega\otimes
\operatorname{id})(W)$, where $\omega=\varphi(\cdot\,a)$.  We will call 
${\mathcal F}(a)$, the {\em Fourier transform\/} of $a$.
\end{defn}

The claim above that ${\mathcal F}(a)\in\hat{M}$ is an easy consequence 
of the fact that $W\in M\otimes\hat{M}$.  Actually, we can be a little 
more precise:

\begin{prop}\label{FT}
Let $a\in\hat{\lambda}(\hat{\mathcal I})$, and let ${\mathcal F}(a)$ be 
the Fourier transform of $a$.  Then we have: ${\mathcal F}(a)\in\lambda
({\mathcal I})\,(\subseteq\hat{M})$.  Moreover, we have:
$$
\bigl\langle\hat{\Lambda}\bigl({\mathcal F}(a)\bigr),\Lambda(x)\bigr\rangle
=\bigl\langle\Lambda(a),\Lambda(x)\bigr\rangle,
$$
for any $x\in{\mathfrak N}_{\varphi}$.  Since the vectors of the form 
$\Lambda(x)$, $x\in{\mathfrak N}_{\varphi}$, are dense in the Hilbert space
${\mathcal H}$, this means that $\hat{\Lambda}\bigl({\mathcal F}(a)\bigr)
=\Lambda(a)$ in ${\mathcal H}$.
\end{prop}

\begin{proof}
As noted in Definition \ref{FTdefn}, we may, at least formally, regard 
${\mathcal F}(a)=(\omega\otimes\operatorname{id})(W)$, where $\omega
=\varphi(\cdot\,a)$.  On the other hand, since $a\in\hat{\lambda}(\hat
{\mathcal I})$, we know that $a\in{\mathfrak N}_{\varphi}$.  So we have: 
$$
\bigl|\omega(x^*)\bigr|=\bigl|\varphi(x^*a)\bigr|=\bigl|\bigl\langle\Lambda(a),
\Lambda(x)\bigr\rangle\bigr|\le L\|\Lambda(x)\|, 
$$
for any $x\in{\mathfrak N}_{\varphi}$ and $L=\|\Lambda(a)\|$.  This means that 
$\omega\in{\mathcal I}$, and we have: ${\mathcal F}(a)=\lambda(\omega)\in
\lambda({\mathcal I})$.  Moreover, from equation \eqref{(phihat)}, we have:
$$
\bigl\langle\hat{\Lambda}\bigl({\mathcal F}(a)\bigr),\Lambda(x)\bigr\rangle
=\omega(x^*)=\varphi(x^*a)=\bigl\langle\Lambda(a),\Lambda(x)\bigr\rangle.
$$
\end{proof}

Let us now define the inverse Fourier transform.  A justification for this 
definition will be given in Theorem~\ref{FIT} below.

\begin{defn}\label{iFTdefn}
For $b\in\lambda({\mathcal I})\,(\subseteq\hat{M})$, define ${\mathcal F}^{-1}(b)
\in M$, such that
$$
{\mathcal F}^{-1}(b):=(\operatorname{id}\otimes\hat{\varphi})
\bigl(W^*(1\otimes b)\bigr).
$$
We will call ${\mathcal F}^{-1}(b)$, the {\em inverse Fourier transform\/} of $b$.
\end{defn}

Formally, it can be written as ${\mathcal F}^{-1}(b)=(\operatorname{id}\otimes
\theta)(W^*)=\hat{\lambda}(\theta)$, where $\theta=\hat{\varphi}(\cdot\,b)$.  So 
by considering equation \eqref{(antipode)}, it may be written (again, formally) 
as ${\mathcal F}^{-1}(b)=S\bigl((\operatorname{id}\otimes\theta)(W)\bigr)
=\bigl(\operatorname{id}\otimes(S\theta)\bigr)(W)$.  In other words, 
${\mathcal F}^{-1}(b)$ may be considered as the linear functional $S\theta
=\theta\bigl(\hat{S}^{-1}(\cdot)\bigr)=\hat{\varphi}\bigl(\hat{S}^{-1}(\cdot)b
\bigr)$.  Compare this with the result in Lemma~2.1 of \cite{VDFourier}, in 
the multiplier Hopf algebra setting.  The reason we have the left Haar weight 
$\hat{\varphi}$ here, instead of the right integral as in \cite{VDFourier}, 
could be attributed to the (``opposite'') way the comultiplication $\hat{\Delta}$ 
is defined in our operator algebra setting: See the remark made in Section~4,
preceding Proposition~\ref{dualpairing}.

\begin{prop}\label{iFT}
Let $b\in\lambda({\mathcal I})$, and let ${\mathcal F}^{-1}(b)$ be the inverse 
Fourier transform of $b$.  Then we have: ${\mathcal F}^{-1}(b)\in\hat{\lambda}
(\hat{\mathcal I})\,(\subseteq M)$.  Moreover, we have:
$$
\bigl\langle\Lambda\bigl({\mathcal F}^{-1}(b)\bigr),\hat{\Lambda}(y)\bigr\rangle
=\bigl\langle\hat{\Lambda}(b),\hat{\Lambda}(y)\bigr\rangle,
$$
for any $y\in{\mathfrak N}_{\hat{\varphi}}$.  Since the vectors of the form 
$\hat{\Lambda}(y)$, $y\in{\mathfrak N}_{\hat{\varphi}}$, is dense in the Hilbert 
space ${\mathcal H}$, this means that $\Lambda\bigl({\mathcal F}^{-1}(b)\bigr)
=\hat{\Lambda}(b)$ in ${\mathcal H}$.
\end{prop}

\begin{proof}
As noted above, we can, at least formally, regard ${\mathcal F}^{-1}(b)
=(\operatorname{id}\otimes\theta)(W^*)$, where $\theta=\hat{\varphi}(\cdot\,b)$. 
The proof that $\theta\in\hat{\mathcal I}$ goes in exactly the same way as before,
so that we have: ${\mathcal F}^{-1}(b)=\hat{\lambda}(\theta)\in\hat{\lambda}
(\hat{\mathcal I})$.  Meanwhile, from equation~\eqref{(phihatdual)}, we have:
$$
\bigl\langle\Lambda\bigl({\mathcal F}^{-1}(b)\bigr),\hat{\Lambda}(y)\bigr\rangle
=\theta(y^*)=\hat{\varphi}(y^*b)=\bigl\langle\hat{\Lambda}(b),\hat{\Lambda}(y)
\bigr\rangle.
$$
\end{proof}

Observe that the definition of ${\mathcal F}^{-1}$ was obtained by imitating 
the definition of ${\mathcal F}$, changing $\varphi$ into $\hat{\varphi}$,
and $\lambda$ into $\hat{\lambda}$.  Still, its proper justification is 
provided by the following result, which says that ${\mathcal F}^{-1}
\bigl({\mathcal F}(a)\bigr)=a$ and ${\mathcal F}\bigl({\mathcal F}^{-1}(b)
\bigr)=b$.  This would be our ``Fourier inversion theorem'':

\begin{theorem}\label{FIT}
Let $(M,\Delta)$ and $(\hat{M},\hat{\Delta})$ be a mutually dual pair 
of locally compact quantum groups, and let $\hat{\lambda}(\hat{\mathcal I})
\subseteq M$ and $\lambda({\mathcal I})\subseteq\hat{M}$ be the (dense) 
subalgebras, on which the Fourier transform, ${\mathcal F}$, and the inverse 
Fourier transform, ${\mathcal F}^{-1}$, are defined.  Then we have: 
\begin{enumerate}
  \item For $a\in\hat{\lambda}(\hat{\mathcal I})$, we have: 
${\mathcal F}^{-1}\bigl({\mathcal F}(a)\bigr)=a$.
  \item For $b\in\lambda({\mathcal I})$, we have: 
${\mathcal F}\bigl({\mathcal F}^{-1}(b)\bigr)=b$.
\end{enumerate}
\end{theorem}

\begin{proof}
From Proposition~\ref{FT}, we know that for all $a\in\hat{\lambda}
(\hat{\mathcal I})$, we have: $\hat{\Lambda}\bigl({\mathcal F}(a)\bigr)
=\Lambda(a)\in{\mathcal H}$.  And, from Proposition~\ref{iFT}, we have, 
for all $b\in\lambda({\mathcal I})$ that: $\Lambda\bigl({\mathcal F}^{-1}
(b)\bigr)=\hat{\Lambda}(b)\in{\mathcal H}$.  Taking $b={\mathcal F}(a)
\in\lambda({\mathcal I})$, and combining the two equations, we obtain:
$$
\Lambda\bigl({\mathcal F}^{-1}({\mathcal F}(a))\bigr)=\hat{\Lambda}
\bigl({\mathcal F}(a)\bigr)=\Lambda(a),
$$
in the Hilbert space ${\mathcal H}$.

Let us now write: $a=\hat{\lambda}(\omega)=(\operatorname{id}\otimes
\omega)(W^*)$ and ${\mathcal F}^{-1}\bigl({\mathcal F}(a)\bigr)
=\hat{\lambda}(\theta)=(\operatorname{id}\otimes\theta)(W^*)$, where 
$\omega, \theta\in\hat{\mathcal I}$.  By equation~\eqref{(phihatdual)},
and by the observation that $\Lambda\bigl({\mathcal F}^{-1}({\mathcal F}(a))
\bigr)=\Lambda(a)$, we see that for any $y\in{\mathfrak N}_{\hat{\varphi}}$,
we have:
$$
\omega(y^*)=\bigl\langle\Lambda(a),\hat{\Lambda}(y)\bigr\rangle
=\bigl\langle\Lambda({\mathcal F}^{-1}({\mathcal F}(a))),\hat{\Lambda}(y)
\bigr\rangle=\theta(y^*).
$$
Since ${\mathfrak N}_{\hat{\varphi}}$ is dense in $\hat{M}$, it follows that 
$\theta=\omega$, in $\hat{M}_*$.  From this, it follows that $\hat{\lambda}
(\theta)=\hat{\lambda}(\omega)$, in $M$.  In other words, ${\mathcal F}^{-1}
\bigl({\mathcal F}(a)\bigr)=a\in\hat{\lambda}(\hat{\mathcal I})\,(\subseteq M)$.

In exactly the same way, we can give the proof for ${\mathcal F}\bigl
({\mathcal F}^{-1}(b)\bigr)=b\in\lambda({\mathcal I})\,(\subseteq\hat{M})$.
\end{proof}

In the below is our version of the ``Plancherel formula'', in terms of the Haar 
weights.  Since the vectors $\Lambda(a)$, $a\in\hat{\lambda}(\hat{\mathcal I})$, 
and $\hat{\Lambda}(b)$, $b\in\lambda({\mathcal I})$, are dense in the Hilbert 
space ${\mathcal H}$, the following proposition implies that the maps 
${\mathcal F}$ and ${\mathcal F}^{-1}$ can be considered as unitary maps 
on ${\mathcal H}$ (as is to be expected).

\begin{prop}\label{plancherel}
Let $\hat{\lambda}(\hat{\mathcal I})\subseteq M$ and $\lambda({\mathcal I})
\subseteq\hat{M}$ be the (dense) subalgebras, as defined earlier.  Then we have: 
\begin{enumerate}
  \item For $a\in\hat{\lambda}(\hat{\mathcal I})\,(\subseteq M)$, we have:
$\hat\varphi\bigl({\mathcal F}(a)^*{\mathcal F}(a)\bigr)=\varphi(a^*a)$.
  \item For $b\in\lambda({\mathcal I})\,(\subseteq\hat{M})$, we have:
$\varphi\bigl({\mathcal F}^{-1}(b)^*{\mathcal F}^{-1}(b)\bigr)=\hat{\varphi}(b^*b)$.
\end{enumerate}
\end{prop}

\begin{proof}
For $a\in\hat{\lambda}(\hat{\mathcal I})$, we saw from Proposition~\ref{FT} that
$\hat{\Lambda}\bigl({\mathcal F}(a)\bigr)=\Lambda(a)\,\in{\mathcal H}$.  It follows 
that:
$$
\hat\varphi\bigl({\mathcal F}(a)^*{\mathcal F}(a)\bigr)=\bigl\langle\hat{\Lambda}
({\mathcal F}(a)),\hat{\Lambda}({\mathcal F}(a))\bigr\rangle=\bigl\langle\Lambda(a),
\Lambda(a)\bigr\rangle=\varphi(a^*a).
$$
The second equation can be proved in exactly the same way.
\end{proof}
  
Meanwhile, remembering from classical analysis that the Fourier transform converts 
multiplication of functions to convolution product, and vice versa, we can use 
our Fourier transform to define the ``convolution product'' on the quantum group 
$(M,\Delta)$:

\begin{defn}\label{convolution}
For $a,c\in\hat{\lambda}(\hat{\mathcal I})\,(\subseteq M)$, define their
{\em convolution product\/}, written $a\ast c$, by
$$
a\ast c:={\mathcal F}^{-1}\bigl({\mathcal F}(a){\mathcal F}(c)\bigr).
$$
\end{defn}

Since ${\mathcal F}(a)$ and ${\mathcal F}(c)$ are contained in $\lambda
({\mathcal I})$, and since $\lambda({\mathcal I})$ is a subalgebra of $\hat{M}$, 
we have ${\mathcal F}(a){\mathcal F}(c)\in\lambda({\mathcal I})$.  So by 
Proposition~\ref{iFT}, the expression ${\mathcal F}^{-1}\bigl({\mathcal F}(a)
{\mathcal F}(c)\bigr)$ is well-defined, and is contained in $\hat{\lambda}
(\hat{\mathcal I})$.  While this is certainly all valid, some more discussion 
will be helpful for further understanding of the convolution product.

\begin{prop}\label{convolutionprod}
Suppose $a,c\in\hat{\lambda}(\hat{\mathcal I})\,\subseteq M$, and let $a\ast c
\in\hat{\lambda}(\hat{\mathcal I})$ be their convolution product.  Then we have 
the following, alternative description:
$$
a\ast c=(\varphi\otimes\operatorname{id})\bigl([(S^{-1}\otimes\operatorname{id})
(\Delta c)](a\otimes1)\bigr).
$$
\end{prop}

\begin{proof}
For $a,c\in\hat{\lambda}(\hat{\mathcal I})$, we know from Proposition~\ref{FT} that 
we can write: ${\mathcal F}(a)=(\omega_1\otimes\operatorname{id})(W)$ and 
${\mathcal F}(c)=(\omega_2\otimes\operatorname{id})(W)$, where $\omega_1
=\varphi(\cdot\,a)\in{\mathcal I}$ and $\omega_2=\varphi(\cdot\,c)\in{\mathcal I}$.
Then by Lemma~\ref{lemsubalg} (see also \cite{BS}), we can write:
$$
{\mathcal F}(a){\mathcal F}(c)=(\omega_1\otimes\operatorname{id})(W)
(\omega_2\otimes\operatorname{id})(W)=(\mu\otimes\operatorname{id})(W)
\,\in\lambda({\mathcal I}),
$$
where $\mu\in M_*$ is given by $\mu(x)=(\omega_1\otimes\omega_2)
\bigl(\Delta(x)\bigr)$.  So for any $x\in{\mathfrak N}_{\varphi}$, we have:
\begin{align}
\mu(x^*)&=(\omega_1\otimes\omega_2)\bigl(\Delta(x^*)\bigr)=(\varphi\otimes
\varphi)\bigl(\Delta(x^*)(a\otimes c)\bigr)  \notag \\
&=\varphi\bigl(\bigl[(\operatorname{id}\otimes\varphi)((\Delta(x^*)(1\otimes c))
\bigl]a\bigr).
\notag
\end{align}
Recall now the ``strong'' left invariance property (Proposition~5.40 of \cite{KuVa}),
which says that for any $x,c\in{\mathfrak N}_{\varphi}$, it is known that 
$(\operatorname{id}\otimes\varphi)\bigl(\Delta(x^*)(1\otimes c)\bigr)$ is contained
in the domain of the antipode $S$, and that
$$
S\bigl((\operatorname{id}\otimes\varphi)(\Delta(x^*)(1\otimes c))\bigr)
=(\operatorname{id}\otimes\varphi)\bigl((1\otimes x^*)(\Delta c)\bigr).
$$
Combining this result together with the result we obtained above, we have,
for $a,c\in\hat{\lambda}(\hat{\mathcal I})$ and any $x\in{\mathfrak N}_{\varphi}$, 
the following:
\begin{align}
\mu(x^*)&=\varphi\bigl(\bigl[(\operatorname{id}\otimes\varphi)((\Delta(x^*)
(1\otimes c))\bigl]a\bigr) 
=\varphi\bigl(\bigl[S^{-1}((\operatorname{id}\otimes\varphi)((1\otimes x^*)
(\Delta c)))\bigr]a\bigr)  \notag \\
&=(\varphi\otimes\varphi)\bigl((1\otimes x^*)[(S^{-1}\otimes\operatorname{id})
(\Delta c)](a\otimes1)\bigr)  \notag \\
&=\varphi\bigl(x^*(\varphi\otimes\operatorname{id})\bigl([(S^{-1}\otimes
\operatorname{id})(\Delta c)](a\otimes1)\bigr)\bigr).    \notag
\end{align}
Letting $z=(\varphi\otimes\operatorname{id})\bigl([(S^{-1}\otimes\operatorname{id})
(\Delta c)](a\otimes1)\bigr)$, we see that: $\mu(x^*)=\varphi(x^*z)$, for any 
$x\in{\mathfrak N}_{\varphi}$.  Since $\mu\in{\mathcal I}$, the Fourier inversion 
theorem implies that in fact, $z\in\hat{\lambda}(\hat{\mathcal I})$ and that 
${\mathcal F}(z)=(\mu\otimes\operatorname{id})(W)$.  Or, ${\mathcal F}(z)
={\mathcal F}(a){\mathcal F}(c)$.

Remembering Definition~\ref{convolution}, we conclude that:
$$
a\ast c={\mathcal F}^{-1}\bigl({\mathcal F}(a){\mathcal F}(c)\bigr)=z
=(\varphi\otimes\operatorname{id})\bigl([(S^{-1}\otimes\operatorname{id})
(\Delta c)](a\otimes1)\bigr).
$$
\end{proof}

With only a slight change in the formulation and the proof, we can define 
the convolution product also on the dual quantum group.  In the below is 
the corresponding result:

\begin{prop}
For $b,d\in\lambda({\mathcal I})\,\subseteq\hat{M}$, we can define their
``convolution product'', written $b\ast d$, by
$$
b\ast d:={\mathcal F}\bigl({\mathcal F}^{-1}(b){\mathcal F}^{-1}(d)\bigr).
$$
Then we have the following, alternative description:
$$
b\ast d=(\hat{\varphi}\otimes\operatorname{id})\bigl([(\hat{S}^{-1}\otimes
\operatorname{id})(\hat{\Delta}d)](b\otimes1)\bigr).
$$
\end{prop}

We skip the proof here, since it is really no different from the one given in 
Proposition~\ref{convolutionprod}.  For this, (3) of Lemma~\ref{lemsubalg}
will be useful.  Meanwhile, see Appendix (Section~5) for a discussion that 
these expressions for the convolution products are natural generalizations 
of the ordinary convolution product in classical analysis.

\section{The dual pairing}

For a finite dimensional Hopf algebra $A$, its dual object is none other than the
dual vector space $A'$, equipped with the Hopf algebra structure obtained naturally
from that of $A$ \cite{Ab}.  In general, however, a typical quantum group would be 
infinite dimensional, and in that case, the dual vector space is too big to be given 
any reasonable structure.  Though there are often ways to get around this problem
to define a dual pairing map, things are more tricky for the analytical settings, 
where the quantum groups are required to have additional ($C^*$-algebra or von 
Neumann algebra) structure.

It turns out that in the locally compact quantum group framework, the dual pairing 
between a quantum group $M$ and its dual $\hat{M}$ is only defined at a dense 
subalgebra level, using the multiplicative unitary operator.  It is certainly 
a correct definition, being a natural generalization of the obvious dual pairing 
between $A$ and $A'$ in the finite-dimensional case.  However, the way the pairing 
is defined makes it rather difficult to work with, and we often have to devise some 
indirect ways to get around this problem.  In this section, we show an alternative 
description of the dual pairing using the Haar weight, motivated by the Fourier 
transform.  This new description of the dual pairing may be useful in some future 
research projects.

Let us begin by recalling the definition of the dual pairing, at the level of 
the dense subalgebras ${\mathcal A}$ and $\hat{\mathcal A}$ (see Section~3) of 
the quantum groups $(M,\Delta)$ and $(\hat{M},\hat{\Delta})$.  That is, for 
$a=(\operatorname{id}\otimes\theta)(W)\in{\mathcal A}$ and  $b=(\omega\otimes
\operatorname{id})(W)\in\hat{\mathcal A}$, we have:
\begin{equation}\label{(pairing)}
\langle b\,|\,a\rangle=\bigl\langle(\omega\otimes\operatorname{id})
(W)\,|\,(\operatorname{id}\otimes\theta)(W)\bigr\rangle:=(\omega\otimes
\theta)(W)=\omega\bigl(a)=\theta(b).
\end{equation}
The definition is suggested by \cite{BS}.  The properties of this dual pairing 
map is given below in Proposition~\ref{dualpairing}.

\begin{rem}
Let us point out here the difference in conventions between pure algebra and 
the operator algebra settings.  In purely algebraic frameworks (Hopf algebras, 
QUE algebras, and even multiplier Hopf algebras), the dual comultiplication on 
$A'$ is simply obtained by dualizing the product on $A$ via the pairing map. 
On the other hand, in the setting of locally compact quantum groups, the 
definition of the dual comultiplication on $\hat{M}$ (as reviewed in Section~2) 
is actually ``flipped''.  This results in the dual pairing given in 
equation~\eqref{(pairing)} to become a ``skew'' pairing, in the sense of 
item~(2) of Proposition~\ref{dualpairing} below, as well as the dual antipode 
in item~(3) below appearing with an inverse.  This ``opposite'' way the 
comultiplication has been chosen on $\hat{M}$ also explains the appearance 
of $\hat{\varphi}$ (as opposed to $\hat{\psi}$) in Definition~\ref{iFTdefn} 
earlier of the inverse Fourier transform.
\end{rem}

\begin{prop}\label{dualpairing}
Let $(M,\Delta)$ and $(\hat{M},\hat{\Delta})$ be the dual pair of locally
compact quantum groups, and let ${\mathcal A}$ and $\hat{\mathcal A}$
be their dense subalgebras, as defined earlier.  Then the map 
$\langle\ \,|\,\ \rangle:\hat{\mathcal A}\times{\mathcal A}\to\mathbb{C}$, 
given by equation \eqref{(pairing)}, is a valid dual pairing.
Moreover, we have:
\begin{enumerate}
 \item $\langle b_1b_2\,|\,a\rangle = \bigl\langle b_1\otimes b_2\,|\,\Delta(a)
\bigr\rangle$, for $a\in{\mathcal A}$, $b_1,b_2\in\hat{\mathcal A}$.
 \item $\langle b\,|\,a_1a_2\rangle = \bigl\langle\hat{\Delta}^{\operatorname
{cop}}(b) \,|\, a_1\otimes a_2\bigr\rangle$, for $a_1,a_2\in{\mathcal A}$, 
$b\in\hat{\mathcal A}$.
 \item $\bigl\langle b\,|\,S(a)\bigr\rangle=\bigl\langle\hat{S}^{-1}(b)\,|\,a\bigr
\rangle$, for $a\in{\mathcal A}$, $b\in\hat{\mathcal A}$.
\end{enumerate}
\end{prop}

\begin{proof}
Bilinearity of $\langle\ ,\ \rangle$ is obvious.  So let us just prove the three
properties.  To prove (1), let $a=(\operatorname{id}\otimes\theta)(W)
\in{\mathcal A}$, and suppose $b_1=(\omega_1\otimes\operatorname{id})(W)
\in\hat{\mathcal A}$ and $b_2=(\omega_2\otimes\operatorname{id})(W)
\in\hat{\mathcal A}$.  Then:
\begin{align}
\langle b_1b_2\,|\,a\rangle
&=\bigl\langle(\omega_1\otimes\operatorname{id})(W)(\omega_2\otimes
\operatorname{id})(W)\,|\,(\operatorname{id}\otimes\theta)(W)\bigr\rangle  
\notag \\ 
&=\bigl\langle(\mu\otimes\operatorname{id})(W)\,|\,(\operatorname{id}\otimes
\theta)(W)\bigr\rangle=(\mu\otimes\theta)(W)  \notag \\
&=\mu\bigl((\operatorname{id}\otimes\theta)(W)\bigr)
=(\omega_1\otimes\omega_2)\bigl(W^*[\operatorname{id}\otimes(\operatorname{id}
\otimes\theta)(W)]W\bigr) \notag \\
&=(\omega_1\otimes\omega_2)(\Delta a)=\bigl\langle b_1\otimes b_2\,|\,\Delta(a)
\bigr\rangle.
\notag
\end{align}
The second equality follows from the multiplicativity of $W$, as shown in 
Lemma~\ref{lemsubalg}.  The definition of $\Delta a$ (for $a\in{\mathcal A}
\subseteq M$) was used in the next to last equality.

For (2), let $a_1=(\operatorname{id}\otimes\theta_1)(W)\in{\mathcal A}$ and 
$a_2=(\operatorname{id}\otimes\theta_2)(W)\in{\mathcal A}$, while $b=(\omega
\otimes\operatorname{id})(W)\in\hat{\mathcal A}$.  Then:
\begin{align}
\langle b\,|\,a_1a_2\rangle
&=\bigl\langle(\omega\otimes\operatorname{id})(W)\,|\,(\operatorname{id}
\otimes\theta_1)(W)(\operatorname{id}\otimes\theta_2)(W)\bigr\rangle  \notag \\ 
&=\bigl\langle(\omega\otimes\operatorname{id})(W)\,|\,(\operatorname{id}
\otimes\nu)(W)\bigr\rangle=(\omega\otimes\nu)(W)=\nu\bigl((\omega\otimes
\operatorname{id})(W)\bigr),
\notag
\end{align}
where $\nu\in\hat{M}_*$ is such that $\nu(y):=(\theta_1\otimes\theta_2)
\bigl(W(y\otimes1)W^*\bigr)$.  This again follows from the multiplicativity of $W$
(see again Lemma~\ref{lemsubalg}).  So we have:
$$
\langle b\,|\,a_1a_2\rangle
=\nu(b)=(\theta_1\otimes\theta_2)\bigl(\hat{\Delta}^{\operatorname{cop}}(b)\bigr)
=\bigl\langle\hat{\Delta}^{\operatorname{cop}}(b)\,|\,a_1\otimes a_2\bigr\rangle.
$$

To prove (3), let $a=(\operatorname{id}\otimes\theta)(W)
\in{\mathcal A}$, and suppose $b=(\omega\otimes\operatorname{id})(W)
\in\hat{\mathcal A}$.  Then by equation~\eqref{(antipodehat)}, we know that 
$(\omega\otimes\operatorname{id})(W^*)$ is contained in ${\mathcal D}(\hat{S})$
and that $\hat{S}\bigl((\omega\otimes\operatorname{id})(W^*)\bigr)=b$.
If, in particular, $\omega\in M_*^{\sharp}$ and $\theta\in\hat{M}_*^{\sharp}$ 
(see Lemma~\ref{lemmasharp}), then we can further write:
$$
\hat{S}^{-1}(b)=(\omega\otimes\operatorname{id})(W^*)=\bigl((\bar{\omega}\otimes
\operatorname{id})(W)\bigr)^*=(\bar{\omega}^{\sharp}\otimes\operatorname{id})(W).
$$
and also: $S(a)=(\operatorname{id}\otimes\theta)(W^*)=(\operatorname{id}\otimes
\bar{\theta^{\sharp}})(W)$.  Then we have:
\begin{align}
\bigl\langle\hat{S}^{-1}(b)\,|\,a\bigr\rangle&=\bigl\langle(\bar{\omega}^{\sharp}
\otimes\operatorname{id}(W)\,|\,(\operatorname{id}\otimes\theta)(W)\bigr\rangle
=(\bar{\omega}^{\sharp}\otimes\theta)(W)  \notag \\
&=\theta\bigl((\bar{\omega}^{\sharp}\otimes\operatorname{id})(W)\bigr)
=\theta\bigl((\omega\otimes\operatorname{id})(W^*)\bigr)
=\omega\bigl((\operatorname{id}\otimes\theta)(W^*)\bigr)   \notag \\
&=\omega\bigl(S[\operatorname{id}\otimes\theta)(W)]\bigr)=\omega\bigl(S(a)\bigr)
=\bigl\langle b\,|\,S(a)\bigr\rangle.
\notag
\end{align}
Since such elements are dense, we conclude that: $\bigl\langle b\,|\,S(a)\bigr
\rangle=\bigl\langle\hat{S}^{-1}(b)\,|\,a\bigr\rangle$, for all $a\in{\mathcal A}$, 
$b\in\hat{\mathcal A}$.
\end{proof}

Except for the appearance of the co-opposite comultiplication 
$\hat{\Delta}^{\operatorname{cop}}$, the proposition shows that 
$\langle\ \,|\,\ \rangle$ is a suitable dual pairing map that corresponds 
to the pairing map on (finite-dimensional) Hopf algebras.  In what follows, 
we will give an alternative description of this pairing map, using Haar 
weights.

\begin{prop}\label{FTpairing}
Let $a\in\hat{\lambda}(\hat{\mathcal I})$ and $b\in\lambda({\mathcal I})$ be
such that we can write ${\mathcal F}(a)=(\omega\otimes\operatorname{id})(W)$,
and ${\mathcal F}^{-1}(b)=(\operatorname{id}\otimes\theta)(W^*)$.  Assume 
further that $\theta\in\hat{M}_*^{\sharp}$ (as in Lemma~\ref{lemmasharp}). 
Then ${\mathcal F}(a)\in\hat{\mathcal A}$, ${\mathcal F}^{-1}(b)\in{\mathcal A}$, 
and we have:
$$
\bigl\langle{\mathcal F}(a)\,|\,{\mathcal F}^{-1}(b)\bigr\rangle
=\bigl\langle\hat{\Lambda}({\mathcal F}(a)),\Lambda({\mathcal F}^{-1}(b)^*)
\bigr\rangle.
$$
\end{prop}

\begin{proof}
It is obvious that ${\mathcal F}(a)=(\omega\otimes\operatorname{id})(W)
\in\hat{\mathcal A}$.  Whereas, $\theta\in\hat{M}_*^{\sharp}$ implies that
${\mathcal F}^{-1}(b)=(\operatorname{id}\otimes\theta)(W^*)
=(\operatorname{id}\otimes\theta^{\sharp})(W^*)^*
=(\operatorname{id}\otimes\bar{\theta^{\sharp}})(W)\in{\mathcal A}$.

Meanwhile, by equation~\eqref{(phihat)} and Proposition~\ref{FT}, we have:
$$
\bigl\langle\hat{\Lambda}({\mathcal F}(a)),\Lambda({\mathcal F}^{-1}(b)^*)
\bigr\rangle=\omega\bigl({\mathcal F}^{-1}(b)\bigr).
$$
In addition, by writing ${\mathcal F}^{-1}(b)^*=(\operatorname{id}\otimes
\theta^{\sharp})(W^*)$, we have:
\begin{align}
\bigl\langle\hat{\Lambda}({\mathcal F}(a)),\Lambda({\mathcal F}^{-1}(b)^*)
\bigr\rangle&=\overline{\bigl\langle\Lambda({\mathcal F}^{-1}(b)^*),\hat{\Lambda}
({\mathcal F}(a))\bigr\rangle}   \notag \\
&=\overline{\bigl\langle\Lambda((\operatorname{id}\otimes\theta^{\sharp})(W^*)),
\hat{\Lambda}({\mathcal F}(a))\bigr\rangle}   \notag \\
&=\overline{\theta^{\sharp}\bigl({\mathcal F}(a)^*\bigr)}
=\bar{\theta^{\sharp}}\bigl({\mathcal F}(a)\bigr),
\notag
\end{align}
where we used equation~\eqref{(phihatdual)} and Proposition \ref{iFT}.
Combining these results, we obtain the following:
\begin{align}
\bigl\langle\hat{\Lambda}({\mathcal F}(a)),\Lambda({\mathcal F}^{-1}(b)^*)
\bigr\rangle&=\omega\bigl({\mathcal F}^{-1}(b)\bigr)
=\bar{\theta^{\sharp}}\bigl({\mathcal F}(a)\bigr)   \notag \\
&=\omega\bigl((\operatorname{id}\otimes\bar{\theta^{\sharp}})(W)\bigr)
=\bar{\theta^{\sharp}}\bigl((\omega\otimes\operatorname{id})(W)\bigr)
\notag \\
&=(\omega\otimes\bar{\theta^{\sharp}})(W)
=\bigl\langle(\omega\otimes\operatorname{id})(W)\,|\,(\operatorname{id}
\otimes\bar{\theta^{\sharp}})(W)\bigr\rangle  \notag \\
&=\bigl\langle{\mathcal F}(a)\,|\,{\mathcal F}^{-1}(b)\bigr\rangle.
\notag
\end{align}
\end{proof}

Observe in the proof above that we needed to work with elements in $\hat{\lambda}
(\hat{M}_*^{\sharp})$ to ensure that taking involution is valid.  As Lemma~\ref
{lemmasharp} shows, this is also closely related to the antipode operation: For 
instance, it can be shown that $\theta^{\sharp}=\bar{\theta}\circ\hat{S}$. 
Considering these, our preferred subspaces from now on will be: $D:=\hat{\lambda}
(\hat{\mathcal I}\cap\hat{M}_*^{\sharp})\,\subseteq{\mathcal A}$, and $\hat{D}:=
\lambda({\mathcal I}\cap M_*^{\sharp})\,\subseteq\hat{\mathcal A}$.  In fact,
$D$ and $\hat{D}$ are dense subalgebras in $M$ and $\hat{M}$, respectively
(See Lemma~\ref{lemmasharp}, and see also Proposition~2.6 of \cite{KuVavN}.).
For elements in $D$ and $\hat{D}$, the following holds:

\begin{cor}
Let $D\subseteq{\mathcal A}$ and $\hat{D}\subseteq\hat{\mathcal A}$ be as 
defined above, and let $a\in D$ and $b\in\hat{D}$.  Then we have the following 
description of the dual pairing, in terms of the inner product:
$$
\langle b\,|\,a\rangle=\bigl\langle\hat{\Lambda}(b),\Lambda(a^*)\bigr\rangle.
$$
\end{cor}

\begin{proof}
By Fourier inversion theorem (see Theorem~\ref{FIT}), we can write:
$a={\mathcal F}^{-1}\bigl({\mathcal F}(a)\bigr)$ and $b={\mathcal F}
\bigl({\mathcal F}^{-1}(b)\bigr)$.  So we are able to use the result of the 
previous proposition.
\end{proof}

Finally, we have the following alternative (up till now not appeared in the 
literatures) description of the dual pairing map, given in terms of the Haar 
weights:

\begin{theorem}\label{pairingnew}
Let $D:=\hat{\lambda}(\hat{\mathcal I}\cap\hat{M}_*^{\sharp})\,\subseteq{\mathcal A}$
and $\hat{D}:=\lambda({\mathcal I}\cap M_*^{\sharp})\,\subseteq\hat{\mathcal A}$ 
be the dense subalgebras defined earlier.  The dual pairing map $\langle\ \,|\,\ \rangle:
\hat{\mathcal A}\times{\mathcal A}\to\mathbb{C}$ given in Proposition~\ref{dualpairing}
takes the following form, if we restrict it to the level of the spaces $D$ and $\hat{D}$:
\begin{align}
\langle b\,|\,a\rangle&=\varphi\bigl(a{\mathcal F}^{-1}(b)\bigr)
=\hat{\varphi}\bigl({\mathcal F}(a^*)^*b\bigr)  \notag \\
&=(\varphi\otimes\hat{\varphi})\bigl[(a\otimes1)W^*(1\otimes b)\bigr].
\notag
\end{align}
Of course, it will satisfy the properties (1), (2), (3) of Proposition~\ref{dualpairing}.
\end{theorem}

\begin{proof}
Recall that the Fourier inversion theorem is very much valid in spaces $D$ and $\hat{D}$.
So for $a\in D$ and $b\in\hat{D}$, we can use the result of Proposition~\ref{FTpairing}
and its Corollary above that $\langle b\,|\,a\rangle=\bigl\langle\hat{\Lambda}(b),
\Lambda(a^*)\bigr\rangle$.  Now note from Proposition~\ref{iFT} that $\hat{\Lambda}(b)
=\Lambda\bigl({\mathcal F}^{-1}(b)\bigr)$ in ${\mathcal H}$.  Thus we have:
$$
\langle b\,|\,a\rangle=\bigl\langle\hat{\Lambda}(b),\Lambda(a^*)\bigr\rangle
=\bigl\langle\Lambda({\mathcal F}^{-1}(b)),\Lambda(a^*)\bigr\rangle
=\varphi\bigl(a{\mathcal F}^{-1}(b)\bigr).
$$
Similarly, we have: $\langle b\,|\,a\rangle=\bigl\langle\hat{\Lambda}(b),\Lambda(a^*)
\bigr\rangle=\bigl\langle\hat{\Lambda}(b),\hat{\Lambda}({\mathcal F}(a^*))\bigr\rangle
=\hat{\varphi}\bigl({\mathcal F}(a^*)^*b\bigr)$.
Since we know from Definition~\ref{iFTdefn} that ${\mathcal F}^{-1}(b)=(\operatorname
{id}\otimes\hat{\varphi})\bigl(W^*(1\otimes b)\bigr)$, and since Definition~\ref{FTdefn} 
implies that ${\mathcal F}(a^*)^*=(\varphi\otimes\operatorname{id})\bigl((a\otimes1)
W^*\bigr)$, we conclude that:
$$
\langle b\,|\,a\rangle
=\varphi\bigl(a{\mathcal F}^{-1}(b)\bigr)
=\hat{\varphi}\bigl({\mathcal F}(a^*)^*b\bigr)
=(\varphi\otimes\hat{\varphi})\bigl[(a\otimes1)W^*(1\otimes b)\bigr].
$$
\end{proof}

\section{Appendix: Case of an ordinary group}

In this Appendix, we will show how some of the results in the earlier sections 
are manifested in the case of an ordinary locally compact group.  Obviously, the 
results here will be mostly familiar.  On the other hand, this exercise will give 
us a clearer understanding of the picture, and will provide a further justification 
of the definitions we chose in the earlier sections.

From now on, let $G$ be a locally compact group, with a fixed left Haar measure, 
$dx$.  Let ${\mathcal H}$ denote the Hilbert space $L^2(G)$.  Let us work with 
the two subalgebras, $M$ and $\hat{M}$, of ${\mathcal B}({\mathcal H})$, as follows.

First consider the commutative von Neumann algebra $M=L^{\infty}(G)$, where 
$a\in L^{\infty}(G)$ is viewed as the multiplication operator $\pi_a$ on 
${\mathcal H}=L^2(G)$, by $\pi_a\xi(x)=a(x)\xi(x)$.  Next consider the group 
von Neumann algebra $\hat{M}={\mathcal L}(G)$, given by the left regular 
representation.  That is, for $b\in C_c(G)$, let $L_b\in{\mathcal B}
({\mathcal H})$ be such that $L_b\xi(x)=\int b(z)\xi(z^{-1}x)\,dz$.  We take 
${\mathcal L}(G)$ to be the $W^*$-closure of $L\bigl(C_c(G)\bigr)$.  These 
are well-known von Neumann algebras, and it is also rather well-known that 
we can give (mutually dual) quantum group structures on them.  We briefly 
review the results below.

Let $W\in{\mathcal B}({\mathcal H}\otimes{\mathcal H})={\mathcal B}\bigl(L^2
(G\times G)\bigr)$ be defined by $W\xi(x,y)=\xi(x,x^{-1}y)$.  It is easy to 
show that it is a multiplicative unitary operator.  We also have:
$$
M=\overline{\bigl\{(\operatorname{id}\otimes\omega)(W):\omega\in{\mathcal B}
({\mathcal H})_*\bigr\}}^w,
$$
and the comultiplication on $M$ is given by $\Delta c=W^*(1\otimes c)W$, 
for $c\in M$.  In effect, this will give us $(\Delta a)(x,y)=a(xy)$, for 
$a\in L^{\infty}(G)$.  The antipode map $S:a\to S(a)$ is such that $S(a)(x)
=a(x^{-1})$, while the left Haar weight is just $\varphi(a)=\int a(x)\,dx$. 
In this way, $(M,\Delta)$ becomes a (commutative) von Neumann algebraic 
quantum group.

Meanwhile, we can show without difficulty that:
$$
\hat{M}=\overline{\bigl\{(\omega\otimes\operatorname{id})(W):\omega
\in{\mathcal B}({\mathcal H})_*\bigr\}}^w,
$$
and the comultiplication on $\hat{M}$ is given by $\hat{\Delta}d=\Sigma 
W(d\otimes1)W^*\Sigma$, for $d\in\hat{M}$.  For $b\in C_c(G)$, this reads: 
$(L\otimes L)_{\hat{\Delta}b}\xi(x,y)=\int b(t)\xi(t^{-1}x,t^{-1}y)\,dt$. 
The antipode map is $\hat{S}:b\to\hat{S}(b)$ such that $\hat{S}(b)(x)
=\delta(x^{-1})b(x^{-1})$, where $\delta$ is the modular function.  And, 
the left Haar weight is given by $\hat{\varphi}(b)=b(1)$, where $1=1_G$ is 
the group identity element.  In this way, we have now a (co-commutative) 
von Neumann algebraic quantum group $(\hat{M},\hat{\Delta})$.

As in Proposition~\ref{dualpairing}, a dual pairing map can be considered 
at the level of certain dense subalgebras of $\hat{M}$ and $M$.  For 
convenience, let us consider $\pi\bigl(C_c(G)\bigr)\subseteq M$ and 
$L\bigl(C_c(G)\bigr)\subseteq\hat{M}$.  The dual pairing defined by 
the multiplicative unitary operator $W$, as given in equation~\eqref
{(pairing)}, becomes:
\begin{equation}\label{(Gpairing)}
\bigl\langle L_b\,|\,\pi_a\bigr\rangle=\int a(x)b(x)\,dx,
\end{equation}
for $\pi_a\in\pi\bigl(C_c(G)\bigr)$ and $b\in L\bigl(C_c(G)\bigr)$.
We will skip the proof here (see Proposition~\ref{pairingG} instead),
though it is really not very difficult, using the operator $W$ defined above.

Now that we are pretty much through with the brief review, let us turn 
our attention to our aim of interpreting the results from earlier sections 
in this current setting of $L^{\infty}(G)$ and ${\mathcal L}(G)$.

\begin{prop}\label{FTiFT}
Let $M=L^{\infty}(G)$ and $\hat{M}={\mathcal L}(G)$ be the mutually dual
quantum groups as above, and consider the subalgebras $\pi\bigl(C_c(G)
\bigr)\subseteq M$ and $L\bigl(C_c(G)\bigr)\subseteq\hat{M}$.  Then the 
Fourier transform and the inverse Fourier transform of Section~3 reads 
as follows:
\begin{enumerate}
  \item For $a\in C_c(G)$, we have: $\pi_a\in M$ and ${\mathcal F}(\pi_a)
=L_a\,\in\hat{M}$.
  \item For $b\in C_c(G)$, we have: $L_b\in\hat{M}$ and ${\mathcal F}^{-1}
(L_b)=\pi_b\,\in M$. 
\end{enumerate}
The Fourier inversion theorem is obvious.
\end{prop}

\begin{proof}
For $a\in C_c(G)$ and any $\xi\in{\mathcal H}$, we have, by Definition~\ref
{FTdefn} that:
$$
{\mathcal F}(\pi_a)\xi(y)=(\varphi\otimes\operatorname{id})\bigl(W(\pi_a
\otimes1)\bigr)\xi(y).
$$
Remembering the definitions of $W$ and $\varphi$ given above, it becomes:
$$
{\mathcal F}(\pi_a)\xi(y)=\int a(x)\xi(x^{-1}y)\,dx=L_a\xi(y).
$$
From this, we obtain: ${\mathcal F}(\pi_a)=L_a$.

Meanwhile, for $b\in C_c(G)$ and any $\xi\in{\mathcal H}$, we have, by 
Definition~\ref{iFTdefn} that:
$$
{\mathcal F}^{-1}(L_b)\xi(x)=(\operatorname{id}\otimes\hat{\varphi})
\bigl(W^*(1\otimes b)\bigr).
$$
Since $W^*\xi(x,y)=\xi(x,xy)$ and since $\hat{\varphi}(b)=b(1)$, this becomes:
$$
{\mathcal F}^{-1}(L_b)\xi(x)=b(x\cdot1)\xi(x)=b(x)\xi(x)=\pi_b\xi(x).
$$
We thus obtain: ${\mathcal F}^{-1}(L_b)=\pi_b$.
\end{proof}

What this proposition shows is that under the Fourier transform, not much 
is happening at the level of functions in $C_c(G)$, and in turn, also at 
the level of the Hilbert space $L^2(G)$ (see Proposition~\ref{planch} below). 
The reason for this ``trivialization'' phenomenon is due to the way the general 
theory of locally compact quantum groups works with the same Hilbert space 
(in this case ${\mathcal H}=L^2(G)$) for both $M$ and $\hat{M}$.  In the 
abelian locally compact group case, this would have been the same as saying 
that $L^2(G)=L^2(\hat{G})$, making the Fourier transform irrelevant.  In a 
sense, this is the reason why the Fourier transform has been ``hidden'' in 
the general theory so far, and its development lacking. 

Next proposition concerns the Plancherel formula, as in Proposition~\ref
{plancherel}.

\begin{prop}\label{planch}
As before, let $a\in C_c(G)$ and $b\in C_c(G)$.  Then:
\begin{enumerate}
  \item $\hat{\varphi}\bigl({\mathcal F}(\pi_a)^*{\mathcal F}(\pi_a)\bigr)
=\varphi(\pi_a^*\pi_a)=\|a\|_2^2$.
  \item $\varphi\bigl({\mathcal F}^{-1}(L_b)^*{\mathcal F}^{-1}(L_b)\bigr)
=\hat{\varphi}(L_b^*L_b)=\|b\|_2^2$.
\end{enumerate}
\end{prop}

\begin{proof}
We just need to use the result of Proposition~\ref{plancherel}.  Meanwhile,
remember that for $a\in L^{\infty}(G)$, the involution is given by $a^*(x)
=\overline{a(x)}$.  Whereas, for $L_b\in{\mathcal L}(G)$, recall: $(L_b)^*
=L_{b^*}$, where $b^*(x)=\delta(x^{-1})\overline{b(x^{-1})}$.
\end{proof}

The convolution product (see Definition~\ref{convolution} and 
Proposition~\ref{convolutionprod}) also takes the familiar form.  The result 
below confirms that the formula given in Proposition~\ref{convolutionprod}
is a natural generalization of the familiar classical convolution product. 
But of course, we could have obtained the same result using 
Definition~\ref{convolution} together with Proposition~\ref{FTiFT}.

\begin{prop}
\begin{enumerate}
  \item For $a,c\in C_c(G)\subseteq L^{\infty}(G)$, their convolution product
reads: $(a\ast c)(y)=\int a(x)c(x^{-1}y)\,dx$.
  \item For $L_b,L_d\in L\bigl(C_c(G))\subseteq{\mathcal L}(G)$, their 
convolution product gives: $(b\ast d)(x)=b(x)d(x)$.
\end{enumerate}
\end{prop}

\begin{proof}
We will use the results from Section~3 and the definitions from the earlier 
part of this Appendix.  Note by the way that in our special case, $S^{-1}=S$ 
and also $\hat{S}^{-1}=\hat{S}$.  So for $a,c\in C_c(G)\subseteq L^{\infty}(G)$,
we have:
$$
(a\ast c)(y)=(\varphi\otimes\operatorname{id})\bigl([(S^{-1}\otimes
\operatorname{id})(\Delta c)](a\otimes1)\bigr)=\int c(x^{-1}y)a(x)\,dx.
$$
And, for $L_b,L_d\in L\bigl(C_c(G))\subseteq{\mathcal L}(G)$ and any
$\xi\in L^2(G)$, observe that:
\begin{align}
L_{b\ast d}\xi(y)&=(\hat{\varphi}\otimes\operatorname{id})\bigl([(\hat{S}^{-1}
\otimes\operatorname{id})(L\otimes L)_{\hat{\Delta}d}](L_b\otimes 1)\bigr)
\xi(y)
\notag \\
&=\int d(t)b(t\cdot1)\xi(t^{-1}y)\,dt=\int b(t)d(t)\xi(t^{-1}y)\,dt=L_{bd}\xi(y),
\notag
\end{align}
where $bd(t)=b(t)d(t)$.
\end{proof}

Finally, let us observe how the new description of the dual pairing given in 
Theorem~\ref{pairingnew} plays out in our special case:

\begin{prop}\label{pairingG}
For $\pi_a\in\pi\bigl(C_c(G)\bigr)\subseteq L^{\infty}(G)$ and $L_b\in L
\bigl(C_c(G)\bigr)\subseteq{\mathcal L}(G)$, the dual pairing can be given
by the formula in Theorem~\ref{pairingnew}:
\begin{align}
\langle L_b\,|\,\pi_a\rangle&=(\varphi\otimes\hat{\varphi})\bigl[(a\otimes1)
W^*(1\otimes b)\bigr]   \notag \\
&=\varphi\bigl(\pi_a{\mathcal F}^{-1}(L_b)\bigr)
=\hat{\varphi}\bigl({\mathcal F}(\pi_a^*)^*L_b\bigr)
=\int a(x)b(x)\,dx,
\notag
\end{align}
which agrees with the equation~\eqref{(Gpairing)}
\end{prop}

\begin{proof}
By Proposition~\ref{FTiFT}, we know that ${\mathcal F}^{-1}(L_b)=\pi_b$.
Therefore, we have: 
$$
\langle L_b\,|\,\pi_a\rangle=\varphi(\pi_a\pi_b)=\varphi(\pi_{ab})
=\int a(x)b(x)\,dx.  
$$
Alternatively, with ${\mathcal F}(\pi_a^*)=L_{\bar{a}}$, we arrive at the 
same result by 
$$
\langle L_b\,|\,\pi_a\rangle=\hat{\varphi}\bigl((L_{\bar{a}})^*L_b\bigr)
=\int\delta(x^{-1})a(x^{-1})b(x^{-1}\cdot1)\,dx=\int a(x)b(x)\,dx.
$$
\end{proof}

As the proof of this last proposition shows, the characterization of the 
dual pairing we obtained in Section~3, via the Haar weight, inner product 
and the Fourier transform (see Corollary of Proposition~\ref{FTpairing}
and Theorem~\ref{pairingnew}), can be quite useful: Compared to the 
approach based only on the definition, the new characterization may often 
provide us a simpler approach to results involving the dual pairing.  This 
aspect may turn out to be useful in the future, in more complicated cases
(for instance, see \cite{BJKwal}).

\bigskip\bigskip\bigskip

\providecommand{\bysame}{\leavevmode\hbox to3em{\hrulefill}\thinspace}
\providecommand{\MR}{\relax\ifhmode\unskip\space\fi MR }
\providecommand{\MRhref}[2]{%
  \href{http://www.ams.org/mathscinet-getitem?mr=#1}{#2}
}
\providecommand{\href}[2]{#2}

\end{document}